\input amstex
\documentstyle{amsppt}
\voffset=-1cm
\magnification=1200
\NoBlackBoxes
\input lpicn
\topmatter
\title  On Homology of Virtual Braids and Burau Representation \endtitle
\author Vladimir V. Vershinin
\endauthor
\address
Sobolev Institute of Mathematics, Novosibirsk, 630090, Russia \endaddress
\email versh\@math.nsc.ru
\endemail
\thanks  This research was supported in part by INTAS Grant 96-0712.
\endthanks
\abstract
Virtual knots arise in the study of Gauss diagrams and Vassiliev invariants
of usual knots. The group of virtual braids on $n$ strings $VB_n$ and its
Burau representation to $GL_n\Bbb Z[t,t^{-1}]$ also can be considered. The
homological properties of the series of groups $VB_n$ and its Burau
representation are studied. The following splitting of infinite loop spaces
is proved for the plus-construction of the classifying space of the virtual
braid group on the infinite number of strings:
$$\Bbb Z\times BVB_\infty^+\simeq\Omega^\infty S^\infty\times S^1
\times Y,$$
where $Y$ is an infinite loop space. Connections with $K_*\Bbb Z$ are
discussed.
\endabstract
\endtopmatter
\document
\specialhead 1. Introduction
\endspecialhead
Virtual knots were introduced recently by L.~Kauffman [K] and studied also by
M.~Goussarov, M.~Polyak and O.~Viro [GPV]. The main motivation lies in the
theory of Gauss diagrams and Gauss codes of knots [K, PV]. Namely, for any
knot diagram it is possible to construct its Gauss diagram and form its Gauss
code. The example of a Gauss diagram of a knot is given on the Figure~1.
\vskip0.5cm

{\def\emline#1#2#3#4#5#6{%
       \put(#1,#2){\special{em:moveto}}%
       \put(#4,#5){\special{em:lineto}}}
\unitlength=.15mm
\special{em:linewidth 3.0pt}
\linethickness{3.0pt}
\latexpic(100., 230.0)(-450, -130)
\emline{-220.}{-40.}{1}{-215.82}{-54.948}{1}
\emline{-215.82}{-54.948}{1}{-211.33}{-67.851}{1}
\emline{-211.33}{-67.851}{1}{-206.57}{-78.798}{1}
\emline{-206.57}{-78.798}{1}{-201.6}{-87.877}{1}
\emline{-201.6}{-87.877}{1}{-196.46}{-95.176}{1}
\emline{-196.46}{-95.176}{1}{-191.2}{-100.78}{1}
\emline{-191.2}{-100.78}{1}{-185.88}{-104.79}{1}
\emline{-185.88}{-104.79}{1}{-180.53}{-107.28}{1}
\emline{-180.53}{-107.28}{1}{-175.22}{-108.35}{1}
\emline{-175.22}{-108.35}{1}{-170.}{-108.07}{1}
\emline{-170.}{-108.07}{1}{-164.91}{-106.55}{1}
\emline{-164.91}{-106.55}{1}{-160.}{-103.87}{1}
\emline{-160.}{-103.87}{1}{-155.32}{-100.12}{1}
\emline{-155.32}{-100.12}{1}{-150.93}{-95.382}{1}
\emline{-150.93}{-95.382}{1}{-146.88}{-89.75}{1}
\emline{-146.88}{-89.75}{1}{-143.2}{-83.311}{1}
\emline{-143.2}{-83.311}{1}{-139.96}{-76.155}{1}
\emline{-139.96}{-76.155}{1}{-137.2}{-68.368}{1}
\emline{-137.2}{-68.368}{1}{-134.97}{-60.04}{1}
\emline{-134.97}{-60.04}{1}{-133.33}{-51.259}{1}
\emline{-133.33}{-51.259}{1}{-132.32}{-42.114}{1}
\emline{-132.32}{-42.114}{1}{-132.}{-32.693}{1}
\emline{-132.}{-32.693}{1}{-132.41}{-23.084}{1}
\emline{-132.41}{-23.084}{1}{-133.6}{-13.376}{1}
\emline{-133.6}{-13.376}{1}{-135.63}{-3.6574}{1}
\emline{-135.63}{-3.6574}{1}{-138.53}{5.9834}{1}
\emline{-138.53}{5.9834}{1}{-142.38}{15.458}{1}
\emline{-142.38}{15.458}{1}{-147.2}{24.678}{1}
\emline{-147.2}{24.678}{1}{-153.06}{33.555}{1}
\emline{-153.06}{33.555}{1}{-160.}{42.}{1}
\emline{-117.86}{-18.971}{1}{-107.}{-7.8815}{1}
\emline{-107.}{-7.8815}{1}{-98.072}{2.46}{1}
\emline{-98.072}{2.46}{1}{-90.971}{12.054}{1}
\emline{-90.971}{12.054}{1}{-85.597}{20.902}{1}
\emline{-85.597}{20.902}{1}{-81.846}{29.005}{1}
\emline{-81.846}{29.005}{1}{-79.618}{36.362}{1}
\emline{-79.618}{36.362}{1}{-78.812}{42.977}{1}
\emline{-78.812}{42.977}{1}{-79.326}{48.848}{1}
\emline{-79.326}{48.848}{1}{-81.057}{53.978}{1}
\emline{-81.057}{53.978}{1}{-83.905}{58.367}{1}
\emline{-83.905}{58.367}{1}{-87.768}{62.016}{1}
\emline{-87.768}{62.016}{1}{-92.544}{64.926}{1}
\emline{-92.544}{64.926}{1}{-98.132}{67.098}{1}
\emline{-98.132}{67.098}{1}{-104.43}{68.533}{1}
\emline{-104.43}{68.533}{1}{-111.34}{69.232}{1}
\emline{-111.34}{69.232}{1}{-118.75}{69.195}{1}
\emline{-118.75}{69.195}{1}{-126.57}{68.424}{1}
\emline{-126.57}{68.424}{1}{-134.69}{66.92}{1}
\emline{-134.69}{66.92}{1}{-143.02}{64.683}{1}
\emline{-143.02}{64.683}{1}{-151.44}{61.714}{1}
\emline{-151.44}{61.714}{1}{-159.87}{58.015}{1}
\emline{-159.87}{58.015}{1}{-168.19}{53.585}{1}
\emline{-168.19}{53.585}{1}{-176.3}{48.427}{1}
\emline{-176.3}{48.427}{1}{-184.12}{42.541}{1}
\emline{-184.12}{42.541}{1}{-191.52}{35.928}{1}
\emline{-191.52}{35.928}{1}{-198.42}{28.589}{1}
\emline{-198.42}{28.589}{1}{-204.7}{20.525}{1}
\emline{-204.7}{20.525}{1}{-210.27}{11.737}{1}
\emline{-210.27}{11.737}{1}{-215.03}{2.2247}{1}
\emline{-215.03}{2.2247}{1}{-218.87}{-8.0096}{1}
\emline{-187.14}{58.971}{1}{-202.17}{62.83}{1}
\emline{-202.17}{62.83}{1}{-215.59}{65.391}{1}
\emline{-215.59}{65.391}{1}{-227.45}{66.744}{1}
\emline{-227.45}{66.744}{1}{-237.8}{66.975}{1}
\emline{-237.8}{66.975}{1}{-246.7}{66.171}{1}
\emline{-246.7}{66.171}{1}{-254.18}{64.422}{1}
\emline{-254.18}{64.422}{1}{-260.31}{61.813}{1}
\emline{-260.31}{61.813}{1}{-265.14}{58.432}{1}
\emline{-265.14}{58.432}{1}{-268.72}{54.368}{1}
\emline{-268.72}{54.368}{1}{-271.09}{49.707}{1}
\emline{-271.09}{49.707}{1}{-272.32}{44.537}{1}
\emline{-272.32}{44.537}{1}{-272.46}{38.946}{1}
\emline{-272.46}{38.946}{1}{-271.54}{33.02}{1}
\emline{-271.54}{33.02}{1}{-269.64}{26.849}{1}
\emline{-269.64}{26.849}{1}{-266.79}{20.518}{1}
\emline{-266.79}{20.518}{1}{-263.05}{14.116}{1}
\emline{-263.05}{14.116}{1}{-258.47}{7.7303}{1}
\emline{-258.47}{7.7303}{1}{-253.11}{1.4482}{1}
\emline{-253.11}{1.4482}{1}{-247.01}{-4.6426}{1}
\emline{-247.01}{-4.6426}{1}{-240.23}{-10.455}{1}
\emline{-240.23}{-10.455}{1}{-232.81}{-15.901}{1}
\emline{-232.81}{-15.901}{1}{-224.81}{-20.893}{1}
\emline{-224.81}{-20.893}{1}{-216.29}{-25.344}{1}
\emline{-216.29}{-25.344}{1}{-207.28}{-29.165}{1}
\emline{-207.28}{-29.165}{1}{-197.85}{-32.271}{1}
\emline{-197.85}{-32.271}{1}{-188.05}{-34.573}{1}
\emline{-188.05}{-34.573}{1}{-177.93}{-35.983}{1}
\emline{-177.93}{-35.983}{1}{-167.53}{-36.414}{1}
\emline{-167.53}{-36.414}{1}{-156.91}{-35.779}{1}
\emline{-156.91}{-35.779}{1}{-146.13}{-33.99}{1}
\special{em:linewidth 1.0pt}
\emline{265.}{-25.}{1}{264.82}{-19.349}{1}
\emline{264.82}{-19.349}{1}{264.29}{-13.72}{1}
\emline{264.29}{-13.72}{1}{263.41}{-8.1357}{1}
\emline{263.41}{-8.1357}{1}{262.17}{-2.6179}{1}
\emline{262.17}{-2.6179}{1}{260.6}{2.8115}{1}
\emline{260.6}{2.8115}{1}{258.68}{8.1312}{1}
\emline{258.68}{8.1312}{1}{256.43}{13.32}{1}
\emline{256.43}{13.32}{1}{253.87}{18.358}{1}
\emline{253.87}{18.358}{1}{250.99}{23.224}{1}
\emline{250.99}{23.224}{1}{247.81}{27.901}{1}
\emline{247.81}{27.901}{1}{244.35}{32.368}{1}
\emline{244.35}{32.368}{1}{240.61}{36.609}{1}
\emline{240.61}{36.609}{1}{236.61}{40.607}{1}
\emline{236.61}{40.607}{1}{232.37}{44.346}{1}
\emline{232.37}{44.346}{1}{227.9}{47.812}{1}
\emline{227.9}{47.812}{1}{223.22}{50.99}{1}
\emline{223.22}{50.99}{1}{218.36}{53.868}{1}
\emline{218.36}{53.868}{1}{213.32}{56.434}{1}
\emline{213.32}{56.434}{1}{208.13}{58.68}{1}
\emline{208.13}{58.68}{1}{202.81}{60.595}{1}
\emline{202.81}{60.595}{1}{197.38}{62.172}{1}
\emline{197.38}{62.172}{1}{191.86}{63.406}{1}
\emline{191.86}{63.406}{1}{186.28}{64.29}{1}
\emline{186.28}{64.29}{1}{180.65}{64.822}{1}
\emline{180.65}{64.822}{1}{175.}{65.}{1}
\emline{175.}{65.}{1}{169.35}{64.822}{1}
\emline{169.35}{64.822}{1}{163.72}{64.29}{1}
\emline{163.72}{64.29}{1}{158.14}{63.406}{1}
\emline{158.14}{63.406}{1}{152.62}{62.172}{1}
\emline{152.62}{62.172}{1}{147.19}{60.595}{1}
\emline{147.19}{60.595}{1}{141.87}{58.68}{1}
\emline{141.87}{58.68}{1}{136.68}{56.434}{1}
\emline{136.68}{56.434}{1}{131.64}{53.868}{1}
\emline{131.64}{53.868}{1}{126.78}{50.99}{1}
\emline{126.78}{50.99}{1}{122.1}{47.812}{1}
\emline{122.1}{47.812}{1}{117.63}{44.346}{1}
\emline{117.63}{44.346}{1}{113.39}{40.607}{1}
\emline{113.39}{40.607}{1}{109.39}{36.609}{1}
\emline{109.39}{36.609}{1}{105.65}{32.368}{1}
\emline{105.65}{32.368}{1}{102.19}{27.901}{1}
\emline{102.19}{27.901}{1}{99.01}{23.224}{1}
\emline{99.01}{23.224}{1}{96.132}{18.358}{1}
\emline{96.132}{18.358}{1}{93.566}{13.32}{1}
\emline{93.566}{13.32}{1}{91.32}{8.1312}{1}
\emline{91.32}{8.1312}{1}{89.405}{2.8115}{1}
\emline{89.405}{2.8115}{1}{87.828}{-2.6179}{1}
\emline{87.828}{-2.6179}{1}{86.594}{-8.1357}{1}
\emline{86.594}{-8.1357}{1}{85.71}{-13.72}{1}
\emline{85.71}{-13.72}{1}{85.178}{-19.349}{1}
\emline{85.178}{-19.349}{1}{85.}{-25.}{1}
\emline{85.}{-25.}{1}{85.178}{-30.651}{1}
\emline{85.178}{-30.651}{1}{85.71}{-36.28}{1}
\emline{85.71}{-36.28}{1}{86.594}{-41.864}{1}
\emline{86.594}{-41.864}{1}{87.828}{-47.382}{1}
\emline{87.828}{-47.382}{1}{89.405}{-52.812}{1}
\emline{89.405}{-52.812}{1}{91.32}{-58.131}{1}
\emline{91.32}{-58.131}{1}{93.566}{-63.32}{1}
\emline{93.566}{-63.32}{1}{96.132}{-68.358}{1}
\emline{96.132}{-68.358}{1}{99.01}{-73.224}{1}
\emline{99.01}{-73.224}{1}{102.19}{-77.901}{1}
\emline{102.19}{-77.901}{1}{105.65}{-82.368}{1}
\emline{105.65}{-82.368}{1}{109.39}{-86.609}{1}
\emline{109.39}{-86.609}{1}{113.39}{-90.607}{1}
\emline{113.39}{-90.607}{1}{117.63}{-94.346}{1}
\emline{117.63}{-94.346}{1}{122.1}{-97.812}{1}
\emline{122.1}{-97.812}{1}{126.78}{-100.99}{1}
\emline{126.78}{-100.99}{1}{131.64}{-103.87}{1}
\emline{131.64}{-103.87}{1}{136.68}{-106.43}{1}
\emline{136.68}{-106.43}{1}{141.87}{-108.68}{1}
\emline{141.87}{-108.68}{1}{147.19}{-110.6}{1}
\emline{147.19}{-110.6}{1}{152.62}{-112.17}{1}
\emline{152.62}{-112.17}{1}{158.14}{-113.41}{1}
\emline{158.14}{-113.41}{1}{163.72}{-114.29}{1}
\emline{163.72}{-114.29}{1}{169.35}{-114.82}{1}
\emline{169.35}{-114.82}{1}{175.}{-115.}{1}
\emline{175.}{-115.}{1}{180.65}{-114.82}{1}
\emline{180.65}{-114.82}{1}{186.28}{-114.29}{1}
\emline{186.28}{-114.29}{1}{191.86}{-113.41}{1}
\emline{191.86}{-113.41}{1}{197.38}{-112.17}{1}
\emline{197.38}{-112.17}{1}{202.81}{-110.6}{1}
\emline{202.81}{-110.6}{1}{208.13}{-108.68}{1}
\emline{208.13}{-108.68}{1}{213.32}{-106.43}{1}
\emline{213.32}{-106.43}{1}{218.36}{-103.87}{1}
\emline{218.36}{-103.87}{1}{223.22}{-100.99}{1}
\emline{223.22}{-100.99}{1}{227.9}{-97.812}{1}
\emline{227.9}{-97.812}{1}{232.37}{-94.346}{1}
\emline{232.37}{-94.346}{1}{236.61}{-90.607}{1}
\emline{236.61}{-90.607}{1}{240.61}{-86.609}{1}
\emline{240.61}{-86.609}{1}{244.35}{-82.368}{1}
\emline{244.35}{-82.368}{1}{247.81}{-77.901}{1}
\emline{247.81}{-77.901}{1}{250.99}{-73.224}{1}
\emline{250.99}{-73.224}{1}{253.87}{-68.358}{1}
\emline{253.87}{-68.358}{1}{256.43}{-63.32}{1}
\emline{256.43}{-63.32}{1}{258.68}{-58.131}{1}
\emline{258.68}{-58.131}{1}{260.6}{-52.812}{1}
\emline{260.6}{-52.812}{1}{262.17}{-47.382}{1}
\emline{262.17}{-47.382}{1}{263.41}{-41.864}{1}
\emline{263.41}{-41.864}{1}{264.29}{-36.28}{1}
\emline{264.29}{-36.28}{1}{264.82}{-30.651}{1}
\emline{264.82}{-30.651}{1}{265.}{-25.}{1}
\emline{252.942}{-70.}{1}{97.0577}{20.}{1}
\emline{97.0577}{-70.}{1}{252.942}{20.}{1}
\emline{175}{-115}{1}{175}{65}{1}
\emline{97.0577}{20.}{1}{109.058}{20.}{1}
\emline{97.0577}{20.}{1}{103.058}{9.6077}{1}
\emline{252.942}{20.}{1}{240.942}{20.}{1}
\emline{252.942}{20.}{1}{246.942}{9.6077}{1}
\emline{175.}{-115.}{1}{181.}{-104.608}{1}
\emline{175.}{-115.}{1}{169.}{-104.608}{1}
\put(215.824, 55.2087){\circle*{10}}
\put(-134.692, 66.9198){\circle*{15}}
\special{em:linewidth 2.0pt}
\emline{-207.581}{16.2213}{1}{-199}{42.}{1}
\emline{-207.581}{16.2213}{1}{-185.5}{31}{1}
\put(-245,-45){{\bf 2}}
\put(-175,70){{\bf 1}}
\put(-120,-45){{\bf 3}}
\put(72,20){{\bf 2}}
\put(165,75){{\bf 1}}
\put(72,-75){{\bf 3}}
\special{em:linewidth 1.0pt}
\emline{-30.}{-25}{1}{30}{-25.}{1}
\emline{30}{-25}{1}{15}{-20.}{1}
\emline{30}{-25}{1}{15}{-30.}{1}
\endlatexpic
}
\vglue.2cm
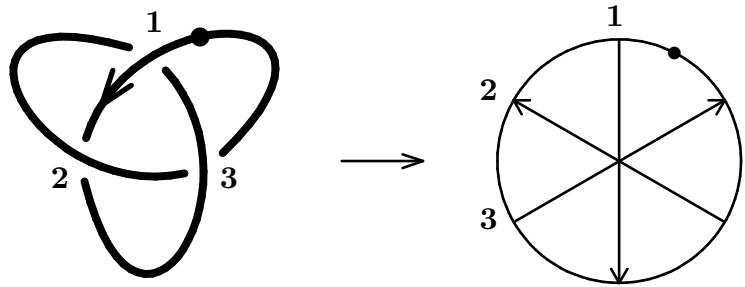
\centerline{ Fig. 1}
\vskip.5cm
\noindent
The Gauss code of this Gauss diagram is the following:

\centerline{ O1U2O3U1O2U3, }
\noindent
where figures mean the crossing points and the symbol ``O" means
overcrossing and the symbol ``U" means undercrossing at this point.
The problem is that not every Gauss diagram (or Gauss code) corresponds
to some knot. To escape this difficulty virtual knots were introduced. Really
virtual knots are those whose Gauss diagram is not a Gauss diagram of any
proper knot. L.~Kauffman generalized many notions from the classical knot
theory to virtual knots, such as fundamental group, rack, quandle, Kauffman
and Jones polynomials. M.~Goussarov, M.~Polyak and O.~Viro [GPV] proved that
the analogues of the upper and the lower presentations of the classical
fundamental group of a knot give two different groups for virtual knots. They
also investigated virtual knots from the point of view of Vassiliev
invariants.

Virtual braids correspond naturally to virtual knots. Some other
generalizations of braids were introduced during the last several years.
These are the braid-permutation group $BP_n$ of R.~Fenn, R.~Rim\'anyi and
C.~Rourke [FRR] and the Baez-Birman monoid $SB_n$ [Ba, Bi2] which embeds into
the singular braid group $SG_n$ [FKR] (see also the description of these
objects in the section 3 of the present papaer). Virtual braids are connected
closely with both of them. They have the generators of the same types and the
relations for virtual braids are of one type of relations for the
braid-permutation group, which do not belong to the Baez-Birman monoid. All
these three types of groups contain the classical braid groups. Burau
representation is defined naturally for the virtual braids and the
braid-permutation group, generalizing the classical Burau representation
[Bu, Bi1].

Cohomologies of classical braid groups appeared in the papers of V.~I.~Arnold
[Arn1, Arn2]. Their study is connected with various mathematical disciplines
and questions. The classical Burau representation was considered from
homotopical point of view by F.~Cohen [C]. Cohomologies of the singular braid
group $SG_n$ and of the braid-permutation group were studied by the author
[Ve]. In the present paper we apply the methods of [Ve] to the investigation
of homological properties of the virtual braid group and its Burau
representation.

\specialhead 2. Virtual Braid Group
\endspecialhead
The following analogues of Reidemeister moves were considered for virtual
knots [K, GPV]:


\centerline{\noindent\latexpic(0,105)(0,-15)
\thicklines
\put(-120,80){\lline(1,-2){20}}
\put(-100,80){\lline(-1,-2){8}}
\put(-70,80){\lline(0,-1){80}}
\put(-50,80){\lline(0,-1){80}}
\put(0,80){\lline(1,-2){40}}
\put(20,80){\lline(-1,-2){8}}
\put(40,80){\lline(-1,-2){18}}
\put(80,80){\lline(1,-2){40}}
\put(100,80){\lline(1,-2){20}}
\put(120,80){\lline(-1,-2){8}}
\put(102,44){\lline(1,2){6}}
\put(-120,40){\lline(1,2){8}}
\put(-120,40){\lline(1,-2){8}}
\put(-100,40){\lline(-1,-2){20}}
\put(0,40){\lline(1,2){8}}
\put(0,40){\lline(1,-2){20}}
\put(120,40){\lline(-1,-2){8}}
\put(12,24){\lline(1,2){6}}
\put(-100,0){\lline(-1,2){8}}
\put(0,0){\lline(1,2){8}}
\put(80,0){\lline(1,2){18}}
\put(100,0){\lline(1,2){8}}
\put(-85,40){\makebox(0,0)[cc]{$\leftrightarrow$}}
\put(60,40){\makebox(0,0)[cc]{$\leftrightarrow$}}
\endlatexpic}
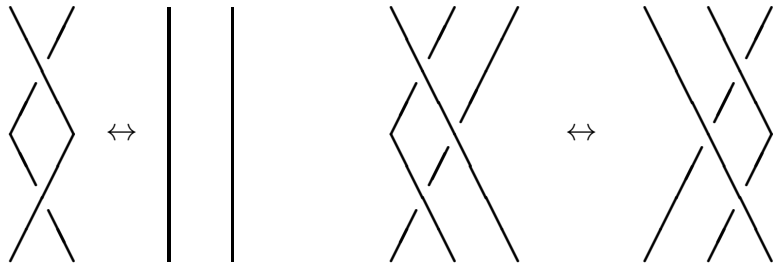
\centerline{Fig. 2}
\vskip 0.2cm

\centerline{\noindent\latexpic(0,110)(0,-15)
\thicklines
\put(-120,80){\lline(1,-2){20}}
\put(-100,80){\lline(-1,-2){20}}
\put(-70,80){\lline(0,-1){80}}
\put(-50,80){\lline(0,-1){80}}
\put(0,80){\lline(1,-2){40}}
\put(20,80){\lline(-1,-2){20}}
\put(40,80){\lline(-1,-2){40}}
\put(80,80){\lline(1,-2){40}}
\put(100,80){\lline(1,-2){20}}
\put(120,80){\lline(-1,-2){40}}
\put(-120,40){\lline(1,-2){20}}
\put(-100,40){\lline(-1,-2){20}}
\put(0,40){\lline(1,-2){20}}
\put(120,40){\lline(-1,-2){20}}
\put(-85,40){\makebox(0,0)[cc]{$\leftrightarrow$}}
\put(60,40){\makebox(0,0)[cc]{$\leftrightarrow$}}
\endlatexpic}
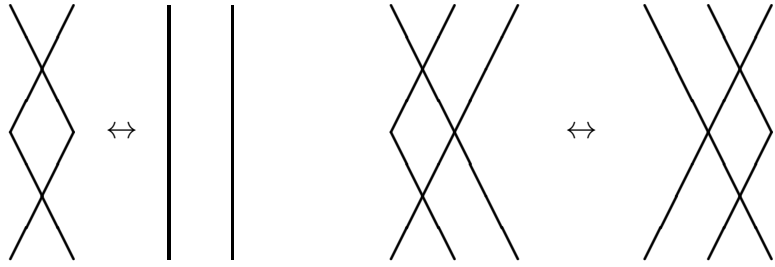
\centerline{Fig. 3}
\vskip 0.2cm

\centerline{\noindent\latexpic(0,105)(-75,-15)
\thicklines
\put(-130,80){\lline(1,-2){40}}
\put(-120,80){\lline(0,-1){80}}
\put(-90,80){\lline(-1,-2){27}}
\put(-60,80){\lline(1,-2){40}}
\put(-30,80){\lline(0,-1){80}}
\put(-20,80){\lline(-1,-2){7}}
\put(-130,0){\lline(1,2){7}}
\put(-60,0){\lline(1,2){27}}
\put(-75,40){\makebox(0,0)[cc]{$\leftrightarrow$}}
\endlatexpic}
\vglue 0.1cm
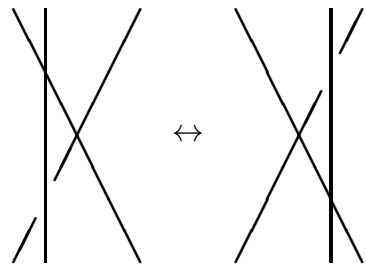
\centerline{Fig. 4}
\vskip 0.2cm

It was emphasized in [GPV] that the following relation does {\it not} fulfill
for virtual knots.

\centerline{\noindent\latexpic(0,110)(75,-15)
\thicklines
\put(20,80){\lline(1,-2){40}}
\put(30,80){\lline(0,-1){13}}
\put(60,80){\lline(-1,-2){17}}
\put(90,80){\lline(1,-2){40}}
\put(120,80){\lline(0,-1){53}}
\put(130,80){\lline(-1,-2){17}}
\put(20,0){\lline(1,2){17}}
\put(30,0){\lline(0,1){53}}
\put(90,0){\lline(1,2){17}}
\put(120,0){\lline(0,1){13}}
\put(75,40){\makebox(0,0)[cc]{$\leftrightarrow$}}
\endlatexpic}
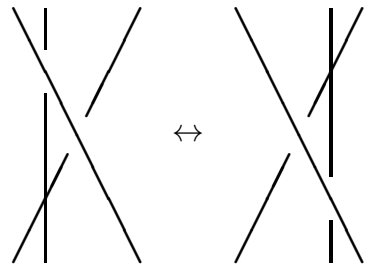
\centerline{Fig. 5. Forbidden move}
\vskip 0.2cm

We can introduce the {\it virtual braid group} by the analogy with the
classical braid group. The difference is that the two types of crossings are
allowed: 1) as usual braids, what is shown at the Figure 6, or 2) as an
intersection of lines on the plane as depicted at the Figure 7.

\centerline{\noindent\latexpic(0,70)(0,-15)
\thicklines
\put(-60,40){\lline(1,-2){20}}
\put(-40,40){\lline(-1,-2){8}}
\put(40,40){\lline(1,-2){8}}
\put(60,40){\lline(-1,-2){20}}
\put(-60,0){\lline(1,2){8}}
\put(60,0){\lline(-1,2){8}}
\endlatexpic}
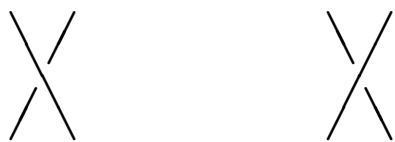
\centerline{Fig. 6}
\vskip 0.3cm

\centerline{\noindent\latexpic(0,50)(0,-15)
\thicklines
\put(-10,40){\lline(1,-2){20}}
\put(10,40){\lline(-1,-2){20}}
\endlatexpic}
\vglue-0.1 cm
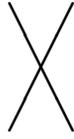
\centerline{Fig. 7}
\vskip 0.5cm

This group is given by the following set of generators:
$\{\zeta_i, \sigma_i, \ \ i=1,2, ..., n-1 \}$ and relations:
$$ \cases \zeta_i^2&=1, \\
\zeta_i \zeta_j &=\zeta_j \zeta_i, \ \ \text {if} \ \ |i-j| >1,\\
\zeta_i \zeta_{i+1} \zeta_i &= \zeta_{i+1} \zeta_i \zeta_{i+1}. \endcases $$
\centerline { The symmetric group relations }

$$ \cases \sigma_i \sigma_j &=\sigma_j \sigma_i, \ \text {if} \  |i-j| >1,
\\ \sigma_i \sigma_{i+1} \sigma_i &= \sigma_{i+1} \sigma_i \sigma_{i+1}.
\endcases \tag1$$
\centerline {The braid group relations }

$$ \cases \sigma_i \zeta_j &=\zeta_j \sigma_i, \ \text {if} \  |i-j| >1,\\
\zeta_i \zeta_{i+1} \sigma_i &= \sigma_{i+1} \zeta_i \zeta_{i+1}.
\endcases $$
\vglue 0.1cm
\centerline {The mixed relations }
\goodbreak

The generator $\sigma_i$ corresponds to the canonical generator of the braid
group $Br_n$ and is depicted at the Figure 8.

\noindent\centerline{\latexpic(0,130)(0,-10)
\thicklines
\put(-100,100){\lline(0,-1){100}}
\put(-50,100){\lline(0,-1){100}}
\put(-25,100){\lline(1,-2){50}}
\put(25,100){\lline(-1,-2){20}}
\put(-25,0){\lline(1,2){20}}
\put(50,100){\lline(0,-1){100}}
\put(100,100){\lline(0,-1){100}}
\put(-100,110){\makebox(0,0)[cc]{$1$}}
\put(-50,110){\makebox(0,0)[cc]{$i-1$}}
\put(-25,110){\makebox(0,0)[cc]{$i$}}
\put(25,110){\makebox(0,0)[cc]{$i+1$}}
\put(50,110){\makebox(0,0)[cc]{$i+2$}}
\put(100,110){\makebox(0,0)[cc]{$n$}}
\put(-75,50){\makebox(0,0)[cc]{.\quad.\quad.}}
\put(75,50){\makebox(0,0)[cc]{.\quad.\quad.}}
\endlatexpic}
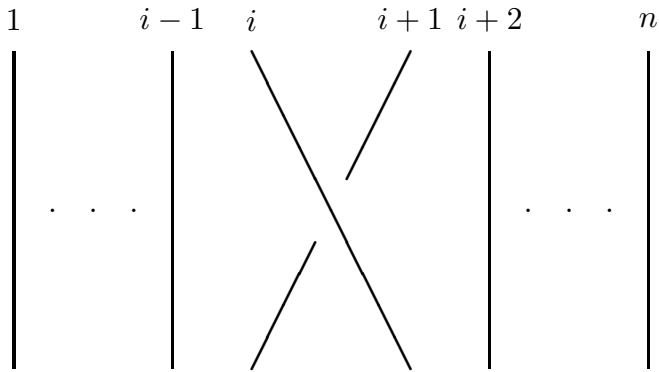
\centerline {Fig. 8}
\vskip 0.5cm

\noindent
The generators $\zeta_i$ correspond to the intersection of lines and is
depicted at the Figure 9.

\noindent\centerline{\latexpic(0,130)(0,-10)
\thicklines
\put(-100,100){\lline(0,-1){100}}
\put(-50,100){\lline(0,-1){100}}
\put(-25,100){\lline(1,-2){50}}
\put(25,100){\lline(-1,-2){50}}
\put(50,100){\lline(0,-1){100}}
\put(100,100){\lline(0,-1){100}}
\put(-100,110){\makebox(0,0)[cc]{$1$}}
\put(-50,110){\makebox(0,0)[cc]{$i-1$}}
\put(-25,110){\makebox(0,0)[cc]{$i$}}
\put(25,110){\makebox(0,0)[cc]{$i+1$}}
\put(50,110){\makebox(0,0)[cc]{$i+2$}}
\put(100,110){\makebox(0,0)[cc]{$n$}}
\put(-75,50){\makebox(0,0)[cc]{.\quad.\quad.}}
\put(75,50){\makebox(0,0)[cc]{.\quad.\quad.}}
\endlatexpic}
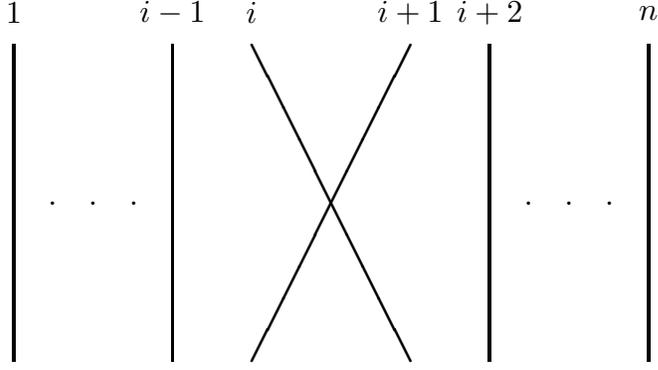
\centerline {Fig. 9}
\vskip 0.3 cm

A homomorphism $j_n$ from the classical braid group $Br_n$
$$j_n: Br_n\to VB_n.$$
is evidently defined by the formulas
$$j_n(\sigma_i)=\sigma_i.$$

We define the category $\Cal V\Cal B$ by analogy with the case of the
classical braids. Its objects $\{\overline 0, \overline  1, ...\} $
correspond to integer numbers from $0$ to infinity
and morphisms are defined by the formula:
$$hom (\overline k,\overline  l)=\cases VB_k, \ \ &\text{if} \quad k=l,\\
\emptyset, \ \  &\text{if} \quad k\not=l. \endcases $$

The pairings
$$\mu_{m,n}: VB_m\times VB_n\to VB_{m+n},$$
are defined in a usual way by the formulas
$$\mu_{m,n}(\sigma_i^\prime)=\sigma_i, \ \mu_{m,n}(\zeta_i^\prime)=\zeta_i;$$
$$\ \sigma_i^\prime,\zeta_i^\prime\in VB_m; \ \sigma_i, \ \zeta_i\in VB_{m+n};$$
$$\mu_{m,n}(\sigma_j^{\prime\prime})=\sigma_{j+m}, \
\mu_{m,n}(\zeta_j^{\prime\prime})=\zeta_{j+m}; $$
$$\ \sigma_i^{\prime\prime}, \ \zeta_i^{\prime\prime}\in VB_n; \
\sigma_{j+m}, \ \zeta_{j+m}\in VB_{m+n}.$$
Geometrically this pairing is defined the same way as for the classical braid
groups: we add (for example on the right) a braid with $n$ strings to a braid
with $m$ strings. It is evident that the pairings for virtual braids agree
with the corresponding pairings for the classical braid groups. They define a
strict monoidal category $\Cal V\Cal B$ and homomorphisms $j_n$ define the
functor from the strict monoidal category $\Cal B$ generated by the braid
groups:
$$\Cal J: \Cal B\to\Cal V\Cal B,$$
which is a morphism of monoidal categories.

We remind that the category $\Cal B$ is a
{\it braided category} as defined by A.~Joyal and R.~Street [JS].
This means that there is given a system of isomorphisms of objects
$$\sigma_{\overline m,\overline n}: \overline{m+ n} \to \overline{n+ m} \ ,$$
natural with respect to morphisms from $\overline m$ to itself and from
$\overline n$ to itself and satisfying the properties of coherence B1 and B2
from [JS]. A braided category becomes symmetric monoidal category if its
braiding satisfies the additional property:
$$\sigma_{\overline m,\overline n}\sigma_{\overline n,\overline m} =
1_{\overline {m+n}} \ . $$
We consider the following system of elements as a braiding $c$ in $\Cal B$
$$\sigma_m ... \sigma_{1}\sigma_{m+1}...\sigma_2...
\sigma_{n+m-1}... \sigma_n \in Br_{m+n}.$$ Graphically it is depicted at
the Figure 10.

\centerline{\noindent\latexpic(0,137)(0,-8)
\thicklines
\put(-50,100){\lline(0,-1){50}}
\put(-25,100){\lline(0,-1){25}}
\put(0,100){\lline(1,-1){50}}
\put(25,100){\lline(-1,-1){10}}
\put(50,75){\lline(-1,-1){10}}
\put(50,100){\lline(0,-1){25}}
\put(-50,50){\lline(1,-1){50}}
\put(-25,75){\lline(1,-1){50}}
\put(0,100){\lline(1,-1){50}}
\put(50,50){\lline(0,-1){50}}
\put(-50,25){\lline(0,-1){25}}
\put(-35,40){\lline(1,1){20}}
\put(25,25){\lline(0,-1){25}}
\put(-10,65){\lline(1,1){19}}
\put(15,40){\lline(1,1){19}}
\put(-50,25){\lline(1,1){10}}
\put(-25,0){\lline(1,1){10}}
\put(-10,15){\lline(1,1){19}}
\put(-25,110){\makebox(0,0)[cc]{$\overbrace{\quad \ \ \ \ \ \ \ \ \ \ \ \ \
\ \ }$}}
\put(-25,120){\makebox(0,0)[cc]{$m$}}
\put(37,110){\makebox(0,0)[cc]{$\overbrace{ \ \ \ \ \ \ \ \ \ \ }$}}
\put(37,120){\makebox(0,0)[cc]{$n$}}
\endlatexpic}
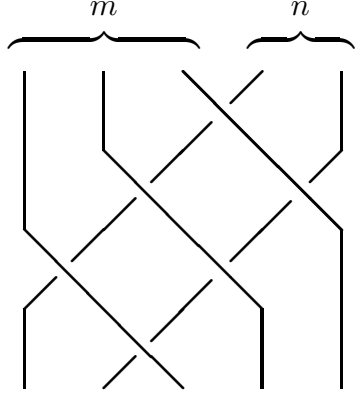
\centerline{Fig. 10}
\vskip 0.3cm
\noindent
Because of the forbidden move (Fig.~5) it does not follow that the image of
the braiding $c$ in the category $\Cal B$ by the functor $\Cal J$ defines a
braiding in the category $\Cal V\Cal B$.

\specialhead 3. The Braid-Permutation Group and the Baez-Birman Monoid
\endspecialhead
Let $F_n$ be the free group of rank $n$ with the set of generators
$\{x_1, ..., x_n \}$ and $\operatorname{Aut} F_n$ be the automorphisms group
of $F_n$. There are the standard inclusions of the symmetric group $\Sigma_n$
and the braid group $Br_n$ into $\operatorname{Aut} F_n$. They can be
described in the following way. Let
$\xi_i \in \operatorname{Aut} F_n, i=1,2, ..., n-1,$ be given by the formula
of the action on the generators:
$$ \cases x_i &\mapsto x_{i+1} \\
x_{i+1} &\mapsto x_i     \\
x_j &\mapsto x_j, j\not=i,i+1. \endcases \tag2$$
Let $ \sigma_i \in \operatorname{Aut} F_n, i=1,2, ..., n-1, $ be given by the
formula of the action on the on generators:
$$ \cases x_i &\mapsto x_{i+1} \\ x_{i+1} &\mapsto x_{i+1}^{-1}x_ix_{i+1} \\
x_j &\mapsto x_j, j\not=i,i+1. \endcases \tag3$$
If we map the standard generators of the symmetric group to $\xi_i$ and the
standard generators of the braid group to $\sigma_i$, then we get the
monomorphisms $t_\Sigma$ and $t_B$ of groups:
$$t_\Sigma: \Sigma_n\rightarrow \operatorname{Aut}F_n,$$
$$t_B: Br_n  \rightarrow \operatorname{Aut}F_n .$$

Let $BP_n$ be the subgroup of $\operatorname{Aut} F_n$, generated by the both
sets of automorphisms $\xi_i$ and $\sigma_i$ of (2) and (3). It is called the
{\it braid-permutation group}. R.~Fenn, R.~Rim\'anyi and C.~Rourke proved
[FRR] that this group is given by the following set of generators:
$\{ \xi_i, \sigma_i, \ \ i=1,2, ..., n-1 \}$ and relations: the braid group
relations for the generators $\sigma_i$ (the same as in (1)), and also:
$$ \cases \xi_i^2&=1, \\
\xi_i \xi_j &=\xi_j \xi_i, \ \ \text {if} \ \ |i-j| >1,\\
\xi_i \xi_{i+1} \xi_i &= \xi_{i+1} \xi_i \xi_{i+1}. \endcases $$
\centerline { The symmetric group relations }

$$ \cases \sigma_i \xi_j &=\xi_j \sigma_i, \ \text {if} \  |i-j| >1,
\\ \xi_i \xi_{i+1} \sigma_i &= \sigma_{i+1} \xi_i \xi_{i+1},
\\ \sigma_i \sigma_{i+1} \xi_i &= \xi_{i+1} \sigma_i \sigma_{i+1}.
\endcases $$
\centerline {The mixed relations }
\medskip

R.~Fenn, R.~Rim\'anyi and C.~Rourke gave the geometrical interpretation of
$BP_n$ as a group of {\it welded braids}.
At first they define a {\it welded braid diagram} on $n$ strings as a
collection of $n$ monotone arcs starting from $n$ points at a horizontal
line of a plane (the top of the diagram) and going down to $n$ points at
another horizontal line (the bottom of the diagram). It is allowed for them
to have crossings of two types: 1) as usual braids, what is shown at the
Figure 6, or 2) to have welds as depicted at the Figure 11.

\centerline{\noindent\latexpic(0,60)(0,-10)
\thicklines
\put(0,20){\circle*{5}}
\put(-10,40){\lline(1,-2){20}}
\put(10,40){\lline(-1,-2){20}}
\endlatexpic}
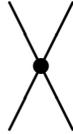
\centerline{Fig. 11}
\vskip 0.5cm
An example of a welded braid diagram is shown at the Figure 12.

\centerline{\noindent\latexpic(0,130)(0,-15)
\thicklines
\put(-17,67){\circle*{5}}
\put(-50,100){\lline(1,-1){100}}
\put(-50,0){\lline(1,2){50}}
\put(0,0){\lline(1,2){14}}
\put(50,100){\lline(-1,-2){30}}
\endlatexpic}
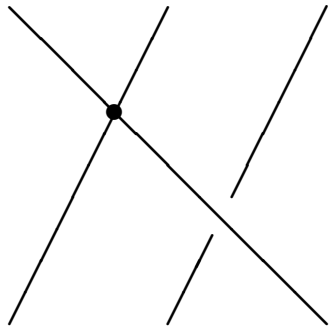
\centerline{Fig. 12}
\vskip 0.5cm

The composition of welded braid diagrams on $n$ strings is defined by
stacking. The diagram with no crossings or welds is an identity with respect
to this composition. So the set of welded braid diagrams on $n$ strings forms
a semi-group denoted by $WD_n$.

R.~Fenn, R.~Rim\'anyi and C.~Rourke define the following types of allowable
moves on welded braid diagrams. They are depicted at the Figures 13-16.

\centerline{\noindent\latexpic(0,120)(0,-15)
\thicklines
\put(-120,80){\lline(1,-2){20}}
\put(-100,80){\lline(-1,-2){8}}
\put(-70,80){\lline(0,-1){80}}
\put(-50,80){\lline(0,-1){80}}
\put(0,80){\lline(1,-2){40}}
\put(20,80){\lline(-1,-2){8}}
\put(40,80){\lline(-1,-2){18}}
\put(80,80){\lline(1,-2){40}}
\put(100,80){\lline(1,-2){20}}
\put(120,80){\lline(-1,-2){8}}
\put(102,44){\lline(1,2){6}}
\put(-120,40){\lline(1,2){8}}
\put(-120,40){\lline(1,-2){8}}
\put(-100,40){\lline(-1,-2){20}}
\put(0,40){\lline(1,2){8}}
\put(0,40){\lline(1,-2){20}}
\put(120,40){\lline(-1,-2){8}}
\put(12,24){\lline(1,2){6}}
\put(-100,0){\lline(-1,2){8}}
\put(0,0){\lline(1,2){8}}
\put(80,0){\lline(1,2){18}}
\put(100,0){\lline(1,2){8}}
\put(-85,40){\makebox(0,0)[cc]{$\leftrightarrow$}}
\put(60,40){\makebox(0,0)[cc]{$\leftrightarrow$}}
\endlatexpic}
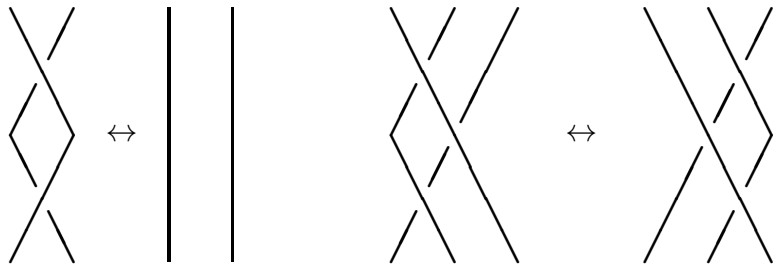
\centerline{Fig. 13}
\vskip 0.5cm

\centerline{\noindent\latexpic(0,110)(0,-15)
\thicklines
\put(-110,60){\circle*{5}}
\put(10,60){\circle*{5}}
\put(110,60){\circle*{5}}
\put(20,40){\circle*{5}}
\put(100,40){\circle*{5}}
\put(-110,20){\circle*{5}}
\put(10,20){\circle*{5}}
\put(110,20){\circle*{5}}
\put(-120,80){\lline(1,-2){20}}
\put(-100,80){\lline(-1,-2){20}}
\put(-70,80){\lline(0,-1){80}}
\put(-50,80){\lline(0,-1){80}}
\put(0,80){\lline(1,-2){40}}
\put(20,80){\lline(-1,-2){20}}
\put(40,80){\lline(-1,-2){40}}
\put(80,80){\lline(1,-2){40}}
\put(100,80){\lline(1,-2){20}}
\put(120,80){\lline(-1,-2){40}}
\put(-120,40){\lline(1,-2){20}}
\put(-100,40){\lline(-1,-2){20}}
\put(0,40){\lline(1,-2){20}}
\put(120,40){\lline(-1,-2){20}}
\put(-85,40){\makebox(0,0)[cc]{$\leftrightarrow$}}
\put(60,40){\makebox(0,0)[cc]{$\leftrightarrow$}}
\endlatexpic}
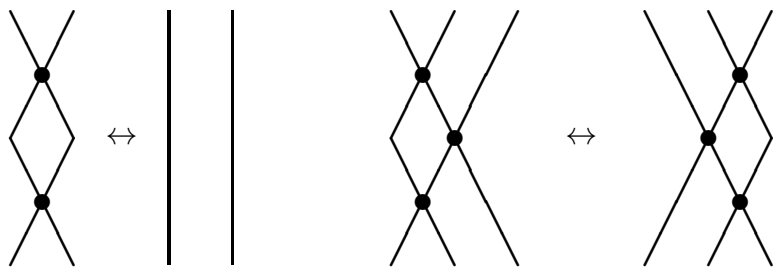
\centerline{Fig. 14}

\centerline{\noindent\latexpic(0,110)(0,-15)
\thicklines
\put(-120,60){\circle*{5}}
\put(120,60){\circle*{5}}
\put(-110,40){\circle*{5}}
\put(-40,40){\circle*{5}}
\put(-30,20){\circle*{5}}
\put(30,20){\circle*{5}}
\put(-130,80){\lline(1,-2){40}}
\put(-120,80){\lline(0,-1){80}}
\put(-90,80){\lline(-1,-2){27}}
\put(-60,80){\lline(1,-2){40}}
\put(-30,80){\lline(0,-1){80}}
\put(-20,80){\lline(-1,-2){7}}
\put(20,80){\lline(1,-2){40}}
\put(30,80){\lline(0,-1){13}}
\put(60,80){\lline(-1,-2){17}}
\put(90,80){\lline(1,-2){40}}
\put(120,80){\lline(0,-1){53}}
\put(130,80){\lline(-1,-2){17}}
\put(-130,0){\lline(1,2){7}}
\put(-60,0){\lline(1,2){27}}
\put(20,0){\lline(1,2){17}}
\put(30,0){\lline(0,1){53}}
\put(90,0){\lline(1,2){17}}
\put(120,0){\lline(0,1){13}}
\put(-75,40){\makebox(0,0)[cc]{$\leftrightarrow$}}
\put(75,40){\makebox(0,0)[cc]{$\leftrightarrow$}}
\endlatexpic}
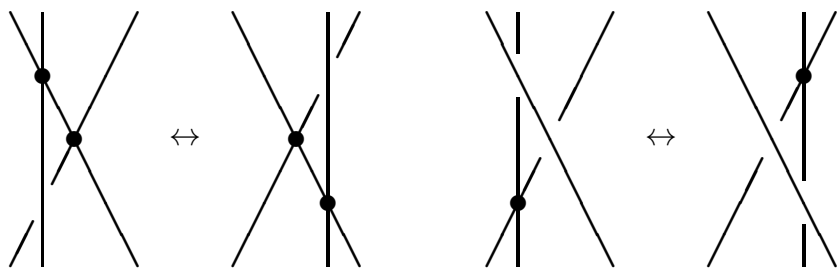
\centerline{Fig. 15}

\centerline{\noindent\latexpic(0,110)(0,-15)
\thicklines
\put(40,60){\circle*{5}}
\put(-90,20){\circle*{5}}
\put(-100,80){\lline(0,-1){40}}
\put(-80,80){\lline(0,-1){40}}
\put(-50,80){\lline(1,-2){20}}
\put(-30,80){\lline(-1,-2){8}}
\put(30,80){\lline(1,-2){20}}
\put(50,80){\lline(-1,-2){20}}
\put(80,80){\lline(0,-1){40}}
\put(100,80){\lline(0,-1){40}}
\put(-100,40){\lline(1,-2){20}}
\put(-80,40){\lline(-1,-2){20}}
\put(-50,40){\lline(1,2){8}}
\put(-50,40){\lline(0,-1){40}}
\put(-30,40){\lline(0,-1){40}}
\put(30,40){\lline(0,-1){40}}
\put(50,40){\lline(0,-1){40}}
\put(80,40){\lline(1,-2){20}}
\put(100,40){\lline(-1,-2){8}}
\put(80,0){\lline(1,2){8}}
\put(0,40){\makebox(0,0)[cc]{$\leftrightarrow$}}
\put(-110,40){\makebox(0,0)[cc]{$. . .$}}
\put(-65,40){\makebox(0,0)[cc]{$. . .$}}
\put(-20,40){\makebox(0,0)[cc]{$. . .$}}
\put(20,40){\makebox(0,0)[cc]{$. . .$}}
\put(65,40){\makebox(0,0)[cc]{$. . .$}}
\put(110,40){\makebox(0,0)[cc]{$. . .$}}
\endlatexpic}
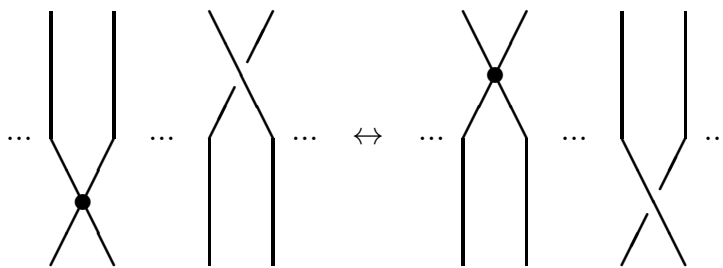
\centerline{Fig. 16}
\vskip 0.5cm

\noindent\centerline{\latexpic(0,130)(0,-10)
\thicklines
\put(0,50){\circle*{5}}
\put(-100,100){\lline(0,-1){100}}
\put(-50,100){\lline(0,-1){100}}
\put(-25,100){\lline(1,-2){50}}
\put(25,100){\lline(-1,-2){50}}
\put(50,100){\lline(0,-1){100}}
\put(100,100){\lline(0,-1){100}}
\put(-100,110){\makebox(0,0)[cc]{$1$}}
\put(-50,110){\makebox(0,0)[cc]{$i-1$}}
\put(-25,110){\makebox(0,0)[cc]{$i$}}
\put(25,110){\makebox(0,0)[cc]{$i+1$}}
\put(50,110){\makebox(0,0)[cc]{$i+2$}}
\put(100,110){\makebox(0,0)[cc]{$n$}}
\put(-75,50){\makebox(0,0)[cc]{.\quad.\quad.}}
\put(75,50){\makebox(0,0)[cc]{.\quad.\quad.}}
\endlatexpic}
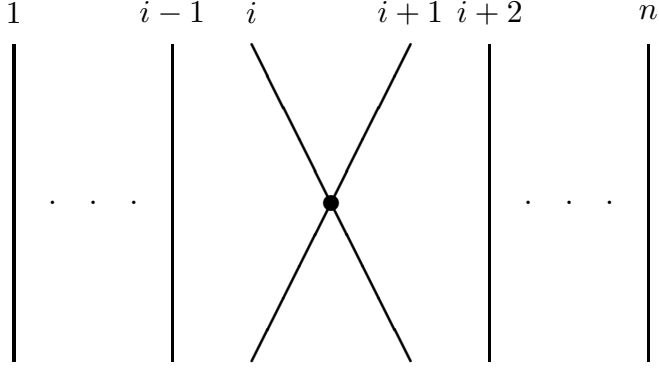
\centerline {Fig. 17}
\vskip 0.3 cm

\noindent
The move from the Figure 16 is the geometric form of the commutativity from
the mixed relations. There are also analogous moves corresponding to the
commutativity from the symmetric group and the braid group relations.

A {\it welded braid} is defined as an equivalence class of welded braid
diagrams under allowable moves. R.~Fenn, R.~Rim\'anyi and C.~Rourke proved
that welded braids form a group and this group is isomorphic to the
braid-permutation group $BP_n$. The generator $\sigma_i$ corresponds to the
canonical generator of the braid group $Br_n$ and is depicted at the
Figure~8. The generators $\xi_i$ correspond to the welded braids depicted at
the Figure~17.

It is also possible to consider welded braids as objects of a
3-dimensional space. We regard them as imbedded in the positive half (with
respect to the third coordinate) of a 3-dimensional space, while the welds
are supposed to belong to the two-dimensional plane with the third coordinate
equal to zero. Hence strings are not allowed to move behind welds.

The {\it Baez-Birman monoid} $SB_n$ (or {\it generalized braid monoid}, or
{\it singular braid monoid}) [Ba, Bi2] is defined as a monoid with
generators
$\sigma_i, \  \sigma_i^{-1}, \ a_i, \ i=1, ..., n-1, $ and relations
$$\align &\sigma_i\sigma_j=\sigma_j\sigma_i, \ \text {if} \  |i-j| >1,\\
&a_ia_j=a_ja_i, \ \text {if} \  |i-j| >1,\\
&a_i\sigma_j=\sigma_ja_i, \ \text {if} \  |i-j| \not=1,\\
&\sigma_i \sigma_{i+1} \sigma_i = \sigma_{i+1} \sigma_i \sigma_{i+1},\\
&\sigma_i \sigma_{i+1} a_i = a_{i+1} \sigma_i \sigma_{i+1},\\
&\sigma_{i+1} \sigma_i a_{i+1} = a_i \sigma_{i+1} \sigma_i,\\
&\sigma_i\sigma_i^{-1}=\sigma_i^{-1}\sigma_i =1.
\endalign $$
In pictures $\sigma_i$ corresponds to canonical generator of the braid group
(right-handed crossing) and $a_i$ represents an intersection of the $i$th
and $(i+1)$st strand, just as in Figures 8 and 17. More detailed geometric
interpretation of the Baez-Birman monoid can be found in the paper of
J.~Birman [Bi2].

A homomorphism $k_n$ from the braid group $Br_n$ is evidently defined:
$$k_n: Br_n\to SB_n.$$
R.~Fenn, E.~Keyman and C.~Rourke proved in [FKR] that the Baez-Birman monoid
embeds in a group $SG_n$ which they call the {\it singular braid group}:
$$SB_n\to SG_n.$$
It means that the elements $a_i$ become invertible and all the relations
of $SB_n$ remain true.

We see that in the two types of mixed relations of the braid-permutation
group one comes from the virtual braid group and another one comes from the
Baez-Birman monoid. Comparing presentations of the virtual braid group and
the braid-permutation group the evident epimorphism can be defined by
the formulas:
$$p(\zeta_i)=\xi_i,$$
$$p(\sigma_i)=\sigma_i.$$
Considering the braid-permutation group as a subgroup of the automorphism
group of the free group $\operatorname{Aut} F_n$ we see that the homomorphism
$p$ defines a representation of the virtual braid group in
$\operatorname{Aut} F_n$.

\specialhead 4. Burau Representation
\endspecialhead
Let us map the generators of the braid-permutation group $BP_n$ to the
following elements of the group $GL_n\Bbb Z[t,t^{-1}]$
$$\sigma_i \mapsto \left( \matrix
E_{i-1} & 0 & 0 & 0 \\
0 & 1-t & t & 0 \\
0 & 1   & 0 & 0 \\
0 & 0   & 0 & E_{n-i-1}
\endmatrix
\right),\tag4$$
$$\xi_i \mapsto \left( \matrix
E_{i-1} & 0 & 0 & 0 \\
0 & 0 & 1 & 0 \\
0 & 1   & 0 & 0 \\
0 & 0   & 0 & E_{n-i-1}
\endmatrix
\right).\tag5$$
\proclaim{Proposition 1} The formulas (4) and (5) define correctly the
representation of the braid-permutation group in $GL_n\Bbb Z[t,t^{-1}]$:
$$r: BP_n\to GL_n\Bbb Z[t,t^{-1}].$$
\endproclaim
\demo{Proof} If follows from the Fox free differential calculus or can be
checked directly.
$\square$
\enddemo
We denote by $BuP_n$ the image in $GL_n\Bbb Z[t,t^{-1}]$ of the
braid-permutation group, which is also the image of the virtual braid group.
So we have the homomorphisms:
$$VB_n\buildrel p\over\longrightarrow BP_n\buildrel r\over
\longrightarrow BuP_n.$$
We need the following simple fact about the group $VB_n$.
\proclaim{Proposition 2} The abelianizations of $VP_n$ and $BuP_n$,
\ $2\leq n\leq\infty$, are both equal to
$\Bbb Z/2\oplus\Bbb Z$:
$$ BV_n/[VB_n, VB_n]\cong\Bbb Z/2\oplus\Bbb Z,$$
$$ BuP_n/[BuP_n, BuP_n]\cong\Bbb Z/2\oplus\Bbb Z.$$
The homomorphisms $p$ and $r$ become isomorphisms after the abelianization.
\endproclaim
\demo{Proof}
We follow the lines of the proof of the analogous fact for the
braid-permutation group [Ve]. Let us add new relations to that of the virtual
braid group $VB_n$:
$$\zeta_i=\zeta_j \ \text{for all} \ i \ \text{and} \ j, $$
$$\sigma_i=\sigma_j \ \text{for all} \ i \ \text{and} \ j, $$
and obtain the epimorphism
$$ab_{VB}: VB_n\to\Bbb Z/2\oplus\Bbb Z.$$
On the other hand if we consider the free product of the symmetric and the
braid groups $\Sigma_n*Br_n$, we obtain from the universality that the
homomorphism of abelianization for the group $\Sigma_n*Br_n$ can be defined
as the composition
$$ab_{\Sigma_n*Br_n}: \Sigma_n*Br_n\buildrel ab_{\Sigma}*ab_{Br}\over
\longrightarrow\Bbb Z/2*\Bbb Z\to\Bbb Z/2\oplus\Bbb Z,$$
where the second map is the canonical epimorphism.
Consider the homomorphism $ab_{\Sigma*Br_n}$ as the composition
$$\Sigma_n*Br_n\to VB_n\buildrel ab_{VB}\over\longrightarrow
\Bbb Z/2\oplus\Bbb Z,$$
where the first map is the canonical epimorphism. Again using the
universality we see that $ab_{VB}$ is the abelianization of $VB_n$. For the
group $VB_n$ the assertion is proved. The determinant of the matrix in (5)
is equal to $-t$ and the determinant of the matrix in (4) is equal to $-1$.
So we have an epimorphism to $\Bbb Z/2\oplus\Bbb Z$. The rest of the proof
follows from the first part.
$\square$
\enddemo
There is an epimorphism
$$\alpha_n: VB_n\to\Bbb Z,$$
which is given by the formulas
$$\zeta_i\mapsto 0 \ \roman {for \ all} \ i,$$
$$\sigma_i\mapsto 1 \ \roman {for \ all} \ i.$$
It follows from the relations, that there exists the epimorphism
$$\phi_n: VB_n \rightarrow \Sigma_n,$$
defined by the formulas:
$$\align \phi_n(\zeta_i)&=\zeta_i,\\
\phi_n(\sigma_i)&=\zeta_i. \endalign $$
Its composition with the canonical inclusion $\nu_n$ of $\Sigma_n$ in
$VB_n$ is equal to the identity map of $\Sigma_n$. These homomorphisms
generate maps of classifying spaces $B\nu_n$ and $B\phi_n$, such that
their composition $$B\Sigma_n \rightarrow BVB_n \rightarrow B\Sigma_n$$ is
equal to identity. We have also the inclusion $j_n$ of the braid group
in $VB_n$, which generates the map of classifying spaces:
$$ B Br_n \rightarrow B BV_n.$$
The composition of $j_n$ and $\phi_n$ gives the
canonical epimorphism: $$\tau_n: Br_n \rightarrow \Sigma_n.$$

We denote by $\Cal Z$ a strict monoidal (tensor) category whose objects
$\{\overline 0, \overline 1, ...\} $ correspond to integer numbers from $0$
to infinity and morphisms are defined by the formula:
$$hom (\overline k,\overline l)=\cases \Bbb Z, \ \  &\text{if} \quad  k=l,\\
\emptyset, \ \  &\text{if} \quad k\not=l. \endcases $$
The product in $\Cal Z$ is defined on objects by the sum of nonnegative
numbers and on morphisms by the sum of integer numbers. This category has
a symmetry which is equal to zero element for all $\overline m$ and
$\overline n$. The homomorphisms $\alpha_n$ induce a morphism of permutative
categories
$$A:\Cal V\Cal B\to\Cal Z$$
and maps of classifying spaces
$$B\alpha_n: B VB_n \rightarrow S^1.$$
We denote by $\gamma$ an inclusion of the group $\Bbb Z$ into $Br_n$:
$$\gamma:\Bbb Z\to Br_n.$$
when the generator of the cyclic group is mapped to one of the generators
$\sigma_i$, say, $\sigma_1$: $\gamma(1)=\sigma_1.$
\proclaim{Theorem 1} There exist maps
$$\psi_V: B\Sigma\times BBr\to BVB^+,$$
and
$$\beta_V:\Omega B(\amalg_{n\geq 0} BVB_n)\to\Omega^\infty S^\infty
\times S^1,$$
such that the map $\psi_V$ becomes a loop map after the group completion,
the map $\beta_V$ is an infinite loop map and it splits by the map
$$\Bbb Z\times B\Sigma^+\times S^1\buildrel Id\times(B\gamma)^+\over
\longrightarrow\Bbb Z\times B\Sigma^+\times BBr^+\buildrel Id\times\psi_V^+
\over\longrightarrow\Bbb Z\times BVB^+.$$
If an infinite loop space $Y$ is a fibre of the map $\beta_V$, then it fits
for the following splitting of the infinite loop spaces:
$$\Omega B(\amalg_{n\geq 0} BVB_n)\simeq\Omega^\infty S^\infty\times S^1
\times Y.$$
The same is true for the groups $BuP_n$: the exist the maps $\psi_P$ and
$\beta_P$ with analogous properties and they commute with the maps $\psi_V$
and $\beta_V$ and corresponding maps for the braid-permutation group with
the help of maps induced by $p$ and $r$.
\endproclaim
\demo{Proof}
We follow the lines of work [Ve]. The category $\Cal V\Cal B$ becomes a
permutative category with the symmetry defined by the elements
$$\zeta_{\overline m,\overline n}=\zeta_m ...\zeta_{1}\zeta_{m+1}...
\zeta_2... \zeta_{n+m-1}...\zeta_n \in VB_{m+n}.$$
It is evident that
$$\zeta^2 =1.$$
By definition, the naturality of the symmetry $\zeta$ means that the
following equality $$\zeta{\overline m,\overline n}\cdot
\mu(b_m^\prime,b_n^{\prime\prime})=
\mu(b_n^{\prime\prime},b_m^\prime)\cdot \zeta{\overline m,\overline n},
\quad b_m^\prime\in VB_{m}, \ \ b_n^{\prime\prime}\in VB_n,$$
is fulfilled. This is equivalent to the expression
$$\zeta{\overline m, \overline n}\cdot\mu(b_m^\prime,b_n^{\prime\prime})
\cdot \zeta{\overline m, \overline n}=\mu(b_n^{\prime\prime},b_m^\prime),$$
which means that the
conjugation by the element $\zeta{\overline m, \overline n}$ transforms the
elements of $VB_m\times VB_n$, canonically lying in $VB_{m+n}$, into the
corresponding elements of $VB_n\times VB_m$. The elements
$\zeta{\overline m, \overline n}$ define a symmetry for the category
formed by symmetric groups, so, for checking the naturality of $\zeta$ in
$\Cal V\Cal B$ it remains to verify the naturality for the generators
$\sigma_i,\ 1\leq i\leq m-1, \ m \leq i \leq m+n$. Let us consider the
corresponding conjugation:
$$\zeta_m ... \zeta_{1}\zeta_{m+1}...\zeta_2... \zeta_{n+m-1}... \zeta_n
\sigma_i \zeta_n ... \zeta_{n+m-1} ...\zeta_{2}\zeta_{m+1}\zeta_1 ...
\zeta_m.$$
If $i>n$, we use the relation
$$\zeta_{i-1}\sigma_i\zeta_{i-1}=\zeta_i\sigma_{i-1}\zeta_i.$$
We have:
$$\zeta_m ... \zeta_{1}\zeta_{m+1}...\zeta_2... \zeta_{n+m-1}... \zeta_n
\sigma_i \zeta_n ... \zeta_{n+m-1}...\zeta_{2}\zeta_{m+1} \zeta_1...\zeta_m
=$$
$$\zeta_m ...\zeta_{1}\zeta_{m+1}...\zeta_2... \zeta_{n+m-1}... \zeta_{i-1}
\sigma_i \zeta_{i-1} ... \zeta_{n+m-1}...\zeta_{2}\zeta_{m+1} \zeta_1
...\zeta_m=$$
$$\zeta_m ...\zeta_{1}\zeta_{m+1}...\zeta_2... \zeta_{n+m-1}... \zeta_{i+1}
\sigma_{i-1} \zeta_i \zeta_{i-1}\zeta_{i-1}... \zeta_{n+m-1}...\zeta_{2}
\zeta_{m+1}\zeta_1...\zeta_m=$$
$$...=\sigma_{i-n}.$$
If $i<n$, we use the relation
$$\zeta_{i}\sigma_{i+1}\zeta_{i}=\zeta_{i+}\sigma_{i}\zeta_{i+1}.$$
We have:
$$\zeta_m ... \zeta_{1}\zeta_{m+1}...\zeta_2... \zeta_{n+m-1}... \zeta_n
\sigma_i \zeta_n ... \zeta_{n+m-1}...\zeta_{2} ... \zeta_{m+1} \zeta_1 ...
\zeta_m=$$
$$\zeta_m ... \zeta_{1}\zeta_{m+1}...\zeta_2... \zeta_i\zeta_{i+m}...
\zeta_{i+2}\zeta_{i+1} \sigma_i \zeta_{i+1} ... \zeta_{i+m}...\zeta_{2}
\zeta_{m+1}\zeta_1...\zeta_m=$$
$$\zeta_m ... \zeta_{1}\zeta_{m+1}...\zeta_2... \zeta_i\zeta_{i+m}...
\zeta_{i+2}\zeta_i\sigma_{i+1} \zeta_i\zeta_{i+2}...\zeta_{i+m}...
\zeta_1 ... \zeta_m=$$
$$\zeta_m ... \zeta_{1}\zeta_{m+1}...\zeta_2...
\zeta_{i+1}\zeta_{i+m}... \zeta_{i+2}\sigma_{i+1}\zeta_{i+2}... \zeta_{i+m}
\zeta_{i+1}...\zeta_{2}...\zeta_{m+1}\zeta_1...\zeta_m=...= $$
$$=\sigma_{i+m}.$$
The conditions of coherence are fulfilled trivially. For the condition B1 of
[JS] we have identically:
$$\zeta_m ... \zeta_{1}\zeta_{m+1}...\zeta_2... \zeta_{n+m-1}...
\zeta_n  \, \cdot \, \zeta_{m+n} ... \zeta_{n+1}\zeta_{m+n+1}...
\zeta_{n+2}... \zeta_{n+m+q-1}...\zeta_{n+q}=$$
$$=\zeta_m ... \zeta_{1}\zeta_{m+1}...\zeta_2... \zeta_{n+m+q-1}...
\zeta_{n+q} .$$
For B2 it is also evident:
$$\zeta_{m+n} ... \zeta_{m+1}\zeta_{m+n+1}...\zeta_{m+2}...
\zeta_{n+m+q-1}...\zeta_{m+q}
\, \cdot \, \zeta_m ... \zeta_{1}\zeta_{m+1}...\zeta_2...
\zeta_{m+q-1}...\zeta_q =$$
$$=\zeta_{m+n} ... \zeta_{1}\zeta_{m+n+1}...\zeta_2... \zeta_{m+n+q-1}...
\zeta_q .$$
So, we proved that the category $\Cal V\Cal B$ is a permutative category.
The morphism $A$ is a morphism of permutative categories and hence
induces the map of infinite loop spaces
$$\Omega B(\amalg_{n\geq 0} BVB_n)\to S^1.$$
Analogously the homomorphisms $\phi_n$ generate the morphism of the
permutative categories
$$\Cal V\Cal B\to \Sigma$$
and the corresponding infinite loop map
$$\Omega B(\amalg_{n\geq 0} BVB_n)\to \Omega^\infty S^\infty.$$
We denote by $\beta_V$ the following composition:
$$\Omega B(\amalg_{n\geq 0} BVB_n)\buildrel\roman{diag}\over\longrightarrow
\Omega B(\amalg_{n\geq 0} BVB_n)\times\Omega B(\amalg_{n\geq 0} BVB_n)\to
\Omega^\infty S^\infty\times S^1.$$
The homomorphisms $\nu$ and $j$ induce the map of products of
classifying spaces:
$$B\Sigma\times BBr \rightarrow BVB\times BVB.$$
The space $BVB^+$ is an $H$-space, so there exists a map
$$\mu: BVB^+\times BVB^+ \rightarrow BVB^+.$$
Let us consider now the composition $f$:
$$B\Sigma\times BBr \buildrel B\nu\times Bj \over
\longrightarrow BVB \times BVB\buildrel q\times q\over\longrightarrow
BVB^+\times BVB^+ \buildrel\mu\over\longrightarrow BVB^+,$$
where $q$ is the canonical map to the plus-construction
$$q: BVB\to BVB^+.$$
From the construction we see that the map
$$B\Sigma \buildrel equiv \over \longrightarrow B\Sigma\times * \buildrel
id\times incl\over\longrightarrow B\Sigma\times BBr\buildrel f\over
\longrightarrow BVB^+$$
is homotopic to $qB\nu$, and the map
$$BBr \buildrel equiv \over \longrightarrow *\times BBr \buildrel incl
\times id\over\longrightarrow B\Sigma\times BBr\buildrel f \over
\rightarrow BVB^+$$
is homotopic to $qBj$. Using the fact that
$\Bbb Z\times BVB^+$ is a loop space, we change a little the map
$$qBj: BBr \to BVB^+,$$
and define the map
$$g: BBr\to BVB^+$$
by the formula
$$g(b)=(qB\tau (b))^{-1}\cdot qBj (b), \ b\in BBr.$$
In this case the composition
$$B^+\phi\cdot g: BBr\to B\Sigma^+$$
is homotopic to zero. We denote by $\psi$ the following composition
$$B\Sigma\times BBr\buildrel B\nu\times g\over \longrightarrow BVB
\times BVB\buildrel q\times q\over\longrightarrow BVB^+\times BVB^+
\buildrel\mu\over\longrightarrow BVB^+$$
and by $\chi$ the composition
$$B\Sigma\times S^1 \buildrel id \times B\gamma\over \longrightarrow
B\Sigma\times BBr \buildrel\psi\over\longrightarrow BVB^+.$$
So, the composition
$$B\Sigma\times BBr\buildrel\psi\over \longrightarrow BVB^+ \buildrel B^+
\phi \times B^+\alpha \over \longrightarrow B\Sigma^+\times S^1 $$
is homotopic to the product $q\times B\alpha$ and the composition
$$B\Sigma\times S^1\buildrel \chi\over\longrightarrow BVB^+ \buildrel
B^+\phi \times B^+\alpha \over \longrightarrow B\Sigma^+\times S^1 $$
is the canonical map from a space to its plus-construction. It gives the
following splitting:
$$B\Sigma^+\times S^1\buildrel\chi^+\over\longrightarrow
BVB^+\buildrel B^+\phi\times B^+\alpha\over\longrightarrow
B\Sigma^+\times S^1. $$
The proof for the groups $BuP_n$ is the same. We notice only that
$K_1\Bbb Z[t,t^{-1}]$ is isomorphic to the group of units of the ring
$\Bbb Z[t,t^{-1}]$ which is isomorphic to $\Bbb Z/2\oplus\Bbb Z$:
$$\Bbb Z[t,t^{-1}]^* \cong \Bbb Z/2\oplus\Bbb Z,$$
and the map is given by determinant:
$$GL(\Bbb Z[t,t^{-1}])\buildrel det\over\longrightarrow \Bbb Z[t,t^{-1}]^* \,
.$$
$\square$
\enddemo
It follows from the Proposition~1 that Burau representation defines a map
$$BVB^+\to BGL(\Bbb Z[t,t^{-1}])^+.$$
Let us consider the homomorphism of homotopy groups that it induces.
$$\pi_*(BVB^+)\to \pi_*(BGL(\Bbb Z[t,t^{-1}])^+.\tag6$$
On the left hand side of (6) we have $\pi_*(\Omega^\infty S^\infty\times
S^1\times Y)\cong \pi_*^S\oplus \pi_*(S^1)\oplus \pi_*(Y)$, where $\pi_*^S$
are the stable homotopy groups of spheres, and on the right hand side
$K_*\Bbb Z[t,t^{-1}]\cong K_*\Bbb Z \oplus K_{*-1}\Bbb Z$.
\proclaim{Proposition 3} The map (6) on the level of $\pi_1$ generates an
epimorphism, which is an isomorphism on the direct summand:
$$\pi_1(\Omega^\infty S^\infty\times S^1) \buildrel\simeq \over
\longrightarrow K_1\Bbb Z\oplus K_{0}\Bbb Z\cong \Bbb Z/2\oplus\Bbb Z$$
For the summand $\pi_*(\Omega^\infty S^\infty)$ the map (6) factors
through $K_*\Bbb Z$ where it is equal to the canonical homomorphism
$$\pi_*^S\to K_*\Bbb Z.\tag7$$
\endproclaim
The proof follows directly from constructions. $\square$

The map (7) as well as the  groups $K_*\Bbb Z$ were studied
in algebraic K-theory (see, for example [Q, We]). The ranks of the
groups $K_*\Bbb Z$ were calculated by A.~Borel [Bo]. Results of
V.~Voevodsky [Vo] allowed to C.~Weibel [W] (with some clarifications of
J.~Rognes and C.~Weibel [RW]) to complete the determination of the 2-torsion
in the groups $K_*\Bbb Z$.

Initially D.~Quillen [Q] considered the commutative diagram
$$\CD B\Sigma^+@>>> BO\\ \downarrow&&\downarrow\\
BGL(\Bbb Z)^+@>>> BGL(\Bbb R)
\endCD$$
and used the results of J.~F.~Adams on J-homomorphism [Ad]. When $i$ is
$8k+1$ or $8k+2$, the 2-torsion part of $K_i\Bbb Z$ is equal to $\Bbb Z/2$
and is generated by the images of Adams' elements $\mu_i\in \pi_*^S$. When
$i=8k+3$ the 2-subgroup of the image of J-homomorphism is $\Bbb Z/8$ and is
contained in the 2-subgroup of $K_i\Bbb Z$ equal to $\Bbb Z/16$. When
$i=8k+7$ the 2-subgroup of the image of J-homomorphism is equal to
$\Bbb Z/w_j$ and is isomorphic to the 2-subgroup of $K_i\Bbb Z$, where
$j=4(k+1)$ and $w_j$ is the largest power of 2 dividing $4j$.

\specialhead  Acknowledgements
\endspecialhead
The author is thankful to Misha Polyak, who explained him the work [GPV]
and to John Rognes for useful information. He is also grateful to Sofia
Lambropoulou and all the organizers of the conference Knots-98 in Delphi,
Greece for the work that they carried out in connection with this very
interesting and successful meeting.
\NoBlackBoxes
 \enddocument